\documentclass[12pt,a4paper,psamsfonts]{amsart}
\usepackage{amssymb,amscd,amsxtra,calc}
\usepackage{cmmib57}
\usepackage[dvips,a4paper,bookmarks,bookmarksopen,
colorlinks=false]{hyperref}

\newtheorem{theorem}{Theorem}[section]
\newtheorem{lemma}[theorem]{Lemma}

\theoremstyle{definition}     % italic or bold etc.

\theoremstyle{remark}

\numberwithin{equation}{section}

\begin{document}

\title[Picard group]
{Picard number of the generic fiber of an abelian fibered
hyperk\"ahler manifold}

\author[Keiji Oguiso]{Keiji Oguiso}
\address{Department of Economy, Keio University, Hiyoshi Hokuku-ku
Yokohama, Japan, and Korea Institute for Advanced Study,
207-43 Cheonryangni-2dong, Dongdaemun-gu, Seoul 130-722, Korea
}
\email{oguiso@hc.cc.keio.ac.jp}


\begin{abstract} We shall show that the Picard number of the generic fiber of an abelian fibered hyperk\"ahler manifold over the projective space is always one. We then give a few applications for the Mordell-Weil group. In particular, by deforming O'Grady's 10-dimensional manifold, we construct an abelian fibered hyperk\"ahler manifold of Mordell-Weil rank 20, which is the maximum possible among all known ones.
\end{abstract}

\maketitle

\section{Introduction}

By a {\it hyperk\"ahler manifold} (HK manifold, for short), we mean a simply connected compact K\"ahler manifold $X$ admitting an everywhere non-degenerate holomorphic $2$-form $\sigma_X$ 
such that $H^{0}(X, \Omega_X^2) = \mathbf C \sigma_X$. Note that $X$
is even dimensional. Put $\dim X = 2n$. We call a surjective morphism $f : X \longrightarrow B$ a {\it fibered HK manifold} if in addition $B$ is normal, projective and $f$ has connected fibers. A fundamental result of Matsushita (\cite{M1}, \cite{M2}) says that each fiber of $f$ is Lagrangian. Especially smooth fibers are complex tori of dimension $n$. They are also projective and hence abelian varieties (see eg. \cite{Og2}). We call $f$ an {\it abelian fibered HK manifold} if in addition $f$ admits a bimeromorphic section, i.e., a subvariety $O \subset X$ such that the induced morphism $f\vert O : O \longrightarrow B$ is bimeromorphic. Then, $X$ is projective (see eg. \cite{Og2}) and we can work in the category of algebraic varieties over $\mathbf C$. In particular, we can speak of the generic fiber $X_{\eta}$ of $f$ in the scheme theoretic sense. This is an abelian variety with origin $O \vert X_{\eta}$, defined over the function field $\mathbf C(B)$. Concerning the base space of a fibered HK manifold, there is a conjecture that it is always the projective space $\mathbf P^n$. Hwang \cite{Hw} has shown that this is true if $X$ is projective and $B$ is smooth.
\par
\vskip 4pt
In this note, we study an abelian fibered HK manifold over the projective space. Our main results are Theorems 1.1, 1.2, 1.3 and 1.4.

\begin{theorem} \label{theorem:AHK}
Let $f : X \longrightarrow \mathbf P^n$ be an abelian fibered HK manifold. Let $K = \mathbf C(\mathbf P^n)$ and let $A_K$ be the generic fiber of $f$. Then, $\rho(A_K) = 1$. Here $\rho(A_K)$ is the Picard number of $A_K$ over $K$.
\end{theorem}

It can happen that $\rho(X_t) \ge 2$ for all smooth closed fiber of $f$ (see 
the first step of the proof of Theorem \ref{theorem:MAX}). The statement is of arithmetical nature. Geometrically, it means that two horizontal divisors
on $X$ are proportional in $NS(X)$ up to vertical divisors. Key of our proof is to compare two different ways to discuss deformations of $f$; Matsushita's way \cite{M3} and Voisin's way \cite{Vo}. The argument is very short, but it is quite mysterious, at least for me, that such considerations yield a strong information about the Picard number of the generic fiber. It may be interesting, if possible, to give a more natural proof based on the monodromy representation around the discriminant locus $D \subset \mathbf P^n$, which is of pure codimension one (\cite{HO}, \cite{Hw}).
\par\vskip 4pt

The set $A_K(K)$ of $K$-rational points, or equivalently, the set of rational sections of $f$, forms an abelian group called the {\it Mordell-Weil group} of $f$. We denote it by ${\rm MW}\, (f)$. When ${\rm MW}\, (f)$ is finitely generated, we call $mw(f) := {\rm rank}\, {\rm MW}\, (f)$ the {\it Mordell-Weil rank}. Mordell-Weil group and Mordell-Weil rank are important invariants of an abelian fibered variety. Combining Theorem 1.1 with \cite{HPW} (or \cite{Kh}, \cite{Og2}), we obtain:

\begin{theorem} \label{theorem:MWG}
Let $f : X \longrightarrow \mathbf P^n$ be as in Theorem \ref{theorem:AHK}.
Let $D = \cup_{i =1}^{k} D_{i}$ be the
irreducible decomposition of the
discriminant divisor $D \subset \mathbf P^n$ and let $m_{i}$
be the number of prime divisors in $X$ lying over $D_{i}$. Then the Mordell-Weil
group ${\rm MW}\, (f)$ is a finitely generated abelian group of rank
$$mw(f) = \rho(X) - 2 -
\sum_{i=1}^{k} (m_{i} -1)\,\, .$$
In particular, $2 \le \rho(X) \le b_2(X) -2$ and $mw(f) \le \rho(X) - 2$.
\end{theorem}

Note that the two inequalities in Theorem \ref{theorem:MWG} are valid even if the base space were singular \cite{Og2}. It is natural to ask if these inequalities are optimal. This is also a problem of arithmetical nature, but again by geometry, we obtain the following:

\begin{theorem} \label{theorem:DMW}
Let $f : X \longrightarrow \mathbf P^n$ be as in Theorem \ref{theorem:AHK}. Assume that $X$ admits at least one (not only rational but also) holomorphic section, i.e., a subvariety $O \subset X$ such that $f \vert O : O \longrightarrow \mathbf P^n$ is an isomorphism. Then, for each integer $\rho$ such that $2 \le \rho \le b_2(X) - 2$, there is an abelian fibered HK manifold $f_{\rho} : X_{\rho} \longrightarrow \mathbf P^n$ which is deformation equivalent to $X$ and which satisfies $\rho(X_{\rho}) = \rho$ and $mw(f_{\rho}) = \rho -2$. Moreover, for each $\rho$, such abelian fibered HK manifolds are dense around $X$ in the Kuranishi space with respect to Euclidean topology.
\end{theorem}

This is a generalization of our earlier result \cite{Og1} about Mordell-Weil groups of Jacobian K3 surfaces (see also \cite{Og2} for an intermediate result).
\par
\vskip 4pt
As an application, we obtain:

\begin{theorem} \label{theorem:MAX}
(1) Let $\rho$ be an integer such that $2 \le \rho \le 21$ and let $n \ge 2$ be an integer. Then, for each such $\rho$ and $n$, there is an abelian fibered HK manifold $f : X \longrightarrow \mathbf P^n$ such that $X$ is deformation equivalent to the Hilbert scheme $S^{[n]}$of $n$-points of a K3 surface $S$, $\rho(X) = \rho$ and $mw(f) = \rho -2$.

(2) Let $\rho$ be an integer such that $2 \le \rho \le 22$. Then, for each such $\rho$, there is an abelian fibered HK manifold $f : X \longrightarrow \mathbf P^5$ such that $X$ is deformation equivalent to O'Grady's $10$-dimensional HK manifold $\mathcal M_4$ \cite{OGr}, $\rho(X) = \rho$ and $mw(f) = \rho -2$. In particular, the case $\rho = 22$ gives the largest record of the Mordell-Weil rank $20$, among all the currently known abelian fibered HK manifolds.
\end{theorem}

Note that $S^{[n]}$ is a HK manifold of $\dim\, S^{[n]} = 2n$ and $b_2(S^{[n]}) = 23$ \cite{Be}, and that $\mathcal M_4$ is a HK manifold with largest known second Betti number, $b_2(\mathcal M_4) = 24$ \cite{Ra2}.

We show Theorems \ref{theorem:AHK}, \ref{theorem:MWG} in Section 2
and Theorems \ref{theorem:DMW}, \ref{theorem:MAX} in Section 3.

\section{Proof of Theorems \ref{theorem:AHK} and \ref{theorem:MWG}}

Let $\mathcal P = \{[\sigma] \in \mathbf P(H^2(X, \mathbf C))\, \vert\, (\sigma, \sigma) = 0\,\, ,\,\, (\sigma, \overline{\sigma}) > 0\}$ be the period domain of $X$. By the local Torelli theorem for HK manifolds (see eg. \cite{GHJ} Proposition 22.11), we can (and will) regard the Kuranishi space $\mathcal K$ of $X$ 
as a small neighbourhood of $[\sigma_X] \in \mathcal P$.

Let $F$ be a smooth closed fiber of $f$
and let $\iota : F \longrightarrow X$ be the inclusion map. As $f$ is fibered over $\mathbf P^n$, the deformations of $X$ that preserve the fibration are of codimension $1$ in $\mathcal K$ by Matsushita \cite{M3} Corollary 1.7 (see also \cite{Sa2}). More precisely, it is a hypersurface defined by $(h, \sigma) = 0$, where $h = f^* H$ is the pullback of the hyperplane class
$H$ of $\mathbf P^n$. Here the fact that the base space is $\mathbf P^n$ is important. As fibers of a fibered HK manifold are Lagrangian \cite{M2}, this deformation is part of deformations of $X$ that keep $F$ Lagrangian. By the easier direction of Voisin's result \cite{Vo}, the latter space is of codimension $r$ in $\mathcal K$, where $r = {\rm rank}\, {\rm Im} (\iota^{*} : H^{2}(X, \mathbf Z) \longrightarrow H^{2}(F, \mathbf Z))$. Thus we have $r = 1$.

Let us show that $\rho(A_K) = 1$. For this, it suffices to show that two elements $L_1$ and $L_2$ of $NS(A_K)$ are proportional. As $A_K$ is projective, we may (and will) assume that both are represented by very ample divisors on $A_K$. We can then write $L_i = [H_i] \vert A_K$ ($i = 1$, $2$). 
Here $H_i$ ($i = 1$, $2$)
are effective Cartier divisors on $X$ and $[H_i] \in NS(X)$ are their classes. 

Let $U \subset \mathbf P^n$ be a dense Zariski open subset of $\mathbf P^n$ such that $f$ is smooth over $U$, the rational section $O$ is holomorphic over $U$ and both $H_i$ are flat over $U$. Let $X_U = f^{-1}(U)$. Then, one can speak of the relative Picard scheme ${\rm Pic}(X_U/U)$ and its identity component ${\rm Pic}^0(X_U/U)$. Note that both ${\rm Pic}(X_U/U)$ and ${\rm Pic}^0(X_U/U)$ satisfy the base change property. (For the Picard scheme, see \cite{Gr} or \cite{FK}, Chapter 9.) For simplicity, we denote the line bundles corresponding to $H_i \vert X_U$ ($i =1,2$) by the same letters. Then $H_i\vert X_U \in {\rm Pic}(X_U/U)(U)$. For a pair of integers $(n, m)$, we set, under the additive notation, 
$$H(n,m) = nH_1 \vert X_U - mH_2 \vert X_U \in {\rm Pic}(X_U/U)(U)\, .$$
By the base change property of ${\rm Pic}(X_U/U)$, the element $H(n, m)$ naturally defines a holomorphic section, say $s(n,m)$, of ${\rm Pic}(X_U/U)$ over 
$U$. Explicitly, $s(n, m)$ is defined by:
$$U \ni t \mapsto H(n, m) \vert F_t \in {\rm Pic}(F_t)\, .$$
Here $t$ is not necessarily a closed point. 
 
Let $F_b$ be the fiber of $f$ over a closed point $b$ of $U$. Then, by $r=1$, 
the two classes 
$[H_i \vert F_b]$ ($i = 1$, $2$) are linearly dependent in $NS(F_b)$. Here we used the fact that $[H_i]$ ($i=1,2$) are elements of $NS(X)$ and $[H_i \vert F_b] = \iota^{*}([H_i])$ ($i = 1$, $2$), where $\iota$ is a natural inclusion map of $F_b$ into $X$. Thus, there is a pair of integers $(n_1, n_2) \not= (0,0)$
such that $[H(n_1, n_2) \vert F_b] = 0$ in $NS(F_b)$. Hence  
$s(n_1, n_2)(b) = H(n_1, n_2) \vert F_b \in {\rm Pic}^0(F_b)$.

As ${\rm Pic}^0(X_U/U)$ is a connected component of ${\rm Pic}(X_U/U)$ and 
$s(n_1, n_2)(b) \in {\rm Pic}^0(F_b)$, it follows that $s(n_1, n_2)(t) \in {\rm Pic}^0(F_t)$ for {\it every} point $t$ of $U$. Now, taking $t = \eta$, we have $H(n_1, n_2) \vert A_K = s(n_1, n_2)(\eta) \in {\rm Pic}^0(A_K)$. Hence $n_1L_1 - n_2L_2 = 0$ in $NS(A_K)$. This completes the proof of Theorem \ref{theorem:AHK}. 
\par
\vskip 4pt

Substituting $\rho(\mathbf P^n) = \rho(A_K) = 1$ into the formula \cite{Og2} Theorem 1.1 (or \cite{HPW} Proposition 2.6 and Lemma 2.4, or \cite{Kh}), we obtain Theorem \ref{theorem:MWG}.

\section{Proof of Theorems \ref{theorem:DMW} and \ref{theorem:MAX}}

Let $f : X \longrightarrow \mathbf P^n$ and $O \simeq \mathbf P^n$ be as in Theorem 1.3 and let $F$ be a general closed fiber of $f$. Then $O$ and $F$ are Lagrangian submanifolds of $X$. Consider the deformations $\mathcal X \longrightarrow \mathcal B \subset \mathcal K$ of $f : X \longrightarrow \mathbf P^n$ which keep $F$ and $O$ Lagrangian. Let $U \subset \mathcal K$ be a small neighbourhood of $[\sigma_X]$.
\par
\vskip 4pt

By \cite{Sa2} (see also \cite{Og2} Proof of Proposition 3.3), the
deformation $\mathcal X \longrightarrow \mathcal B$ forms a projective family of abelian fibered HK manifolds $\tilde{f} : \mathcal X \longrightarrow \mathbf P^n \times \mathcal B$, over $\mathcal B$, with holomorphic section $\mathcal O \simeq \mathbf P^n \times \mathcal B$. The base space $\mathcal B$ is a complex submanifold of codimension $2$ in $\mathcal K$, of the form $(a,\sigma) = (s, \sigma) = 0$, where $a$ and $s$
are linearly independent elements in $H^2(X, \mathbf Q)$. For this, we need the harder direction of \cite{Vo} and the fact that $O \simeq \mathbf P^n$.
\par
\vskip 4pt

Let $\tilde{f}_b : \mathcal X_b \longrightarrow \mathbf P^n$ ($b \in \mathcal B$) be the fibration induced by $\tilde{f}$. By the shape of equations defining $\mathcal B$ and by the Lefschetz $(1,1)$-Theorem, we see that $\rho(\mathcal X_b) \ge 2$ for all $b \in \mathcal B$ and $\rho(\mathcal X_b) =2$ for generic $b \in \mathcal B$. Here, the word ``generic" means that the condition is satisfied outside a countable union of hypersurfaces. In fact, $\rho(\mathcal X_b) \ge \rho$ if and only if $\dim \{\ell \in H^{2}(X, \mathbf Q) \vert (\ell, [\omega_{\mathcal X_b}]) = 0\} \ge \rho$, i.e., there are at least $\rho$ linearly independent rational hyperplanes passing through the point $[\sigma_{\mathcal X_b}]$ corresponding to $\mathcal X_b$. Thus, we have in fact a much stronger statement; in any open neighbourhood of $[\sigma_X] \in \mathcal B$, the subset $\{b\,\vert\, \rho(\mathcal X_b) = \rho\}$ is also dense with respect to Euclidean topology, for each integer $\rho$ such that $2 \le \rho \le b_2(X) -2$. (See \cite{Og1} Theorem 1.1 and its proof for more details).
\par
\vskip 4pt

Let $\mathcal D \subset \mathbf P^n \times \mathcal B$ be the discriminant locus of $\tilde{f}$. The discriminant locus of $\tilde{f}_b$ ($b \in \mathcal B$) is then $\mathcal D_b = \mathcal D \cap (\mathbf P^n \times \{b\})$. This is of pure codimension $1$ in $\mathbf P^n$ (\cite{HO} Proposition 3.1 (2), \cite{Hw} Proposition 4.1). Let $k(b)$ be the number of the irreducible components of $\mathcal D_b$ in $\mathbf P^n$. Then $k(b)$ is upper-semicontinuous on $\mathcal B$ with respect to Zariski topology. Choose $p \in \mathcal B \cap U$ such that $k(p)$ is minimal, say $k$. There is then an open neighbourhood $p \in V \subset \mathcal B \cap U$ such that the number of irreducible components of $\mathcal D_b$ ($b \in V$) is $k$, in particular, constant. For the argument from now, we may (and will) assume that the irreducible decomposition of $\mathcal D$ is $\mathcal D = \cup_{i=1}^{k}
\mathcal D_{i}$ and $\mathcal D_b = \cup_{i=1}^{k} \mathcal D_{i, b}$
is the irreducible decomposition of $\mathcal D_b$ ($b \in V$).
\par
\vskip 4pt

Take one general reference point $o \in V$ such that $\rho(\mathcal X_o) = 2$. Then, by Theorem \ref{theorem:MWG}, $\tilde{f}_o^{-1}(\mathcal D_{i, o})$
are all
irreducible. Moreover, they are reduced as $\tilde{f}$ admits a holomorphic section $\mathcal O$. Take then a resolution of singularities $\tilde{\mathcal D}_i$ of $\tilde{f}^{-1}(\mathcal D_{i})$ and denote by $g : \tilde{\mathcal D}_i \longrightarrow V$ the natural morphism. Recall that smoothness of fibers is a non-empty open condition with respect to Zariski topology for a proper morphism from a smooth variety (in characteristic $0$). Then, by the generality of $o$ and by the fact that $\tilde{f}_o^{-1}(\mathcal D_{i, o})$ is irreducible and reduced, we see that $g^{-1}(o)$ is also irreducible and smooth, and the same holds over some Zariski open neighbourhood $W$ of $o$. Thus $\tilde{f}_{b}^{-1}(\mathcal D_{i, b})$ are irreducible and reduced for all $b \in W$. Thus, from Theorem \ref{theorem:MWG}, we have $mw(\tilde{f}_b) = \rho(\mathcal X_b) - 2$ for all $b \in W$. On the other hand, as we already remarked, the set $W_{\rho} = \{b \in W \vert \rho(\mathcal X_b) = \rho\}$ is dense in $W$ (with respect to Euclidean topology) for each $\rho$ with $2 \le \rho \le b_2(X) -2$. Hence, for each such $\rho$, there is $b \in W \subset U$ such that $\rho(\mathcal X_b) = \rho$
and $mw(\tilde{f}_b) = \rho -2$. This completes the proof of Theorem \ref{theorem:DMW}.
\par
\vskip 4pt

Let us show Theorem \ref{theorem:MAX}. Let $f : S \longrightarrow \mathbf P^1$ be an elliptic K3 surface with a section. Then $f$ induces an abelian fibration $\varphi : S^{[n]} \longrightarrow \mathbf P^n$ with holomorphic section. Thus, Theorem \ref{theorem:MAX} (1) follows from Theorem \ref{theorem:DMW}. Note that $\rho(F) \ge 2$ for each smooth close fiber $F$ of $\varphi$, as $F$ is the product of elliptic curves
(cf. Remark after Theorem \ref{theorem:AHK}).
\par
\vskip 4pt

For Theorem \ref{theorem:MAX} (2), it suffices to find one abelian HK manifold $f : X \longrightarrow \mathbf P^5$ such that $X$ is deformation equivalent to O'Grady's $10$-dimensional HK manifold ${\mathcal M}_4$ and $f$ admits a holomorphic section.

Let $S$ be a generic algebraic K3 surface of degree $2$. Then, ${\rm Pic} S = \mathbf Z H$, $H$ is ample, and $(H^2) = 2$. We denote by $\pi : S \longrightarrow \mathbf P^2$ the finite double cover induced by $\vert H \vert$. Let $\overline{Y} := M((0, [2H], 2))$ be the moduli space of semi-stable sheaves with Mukai vector
$(0, [2H], 2)$ on $S$. The space $\overline{Y}$ is singular, but it admits a symplectic resolution $Y$ \cite{LS}. This $Y$ is a HK manifold and is birational to $\mathcal M_4$ (\cite{OGr} Section 4). Hence $Y$ is deformation equivalent to $\mathcal M_4$ by a fundamental result of Huybrechts \cite{Hu}. For each ${\it I} \in M((0, [2H], 2))$, the fitting ideal of ${\it I}$ is the ideal sheaf of some $C \in \vert 2H \vert$ and vice versa (\cite{Ra1} Remark 2.1.1. See also \cite{Sa1}). Thus, we have a natural surjective morphism $\overline{g} : \overline{Y} \longrightarrow \mathbf P^5 \simeq \vert 2H \vert$. Let $U
\subset \vert 2H \vert$ be the open subset consisting of smooth members. If $[C] \in U$, then $C$ is a smooth curve of genus $5$ and the fiber $\overline{g}^{-1}([C])$ is isomorphic ${\rm Pic}^6(C)$, which is an abelian variety of dimension $5$. The same is true for the induced fibration $g : Y \longrightarrow \mathbf P^5
\simeq \vert 2H \vert$. As $C$ admits a unique $g_2^1(C)$ corresponding to the double cover $\pi \vert C$, the map defined by $C \mapsto 3g_2^1(C) \in {\rm Pic}^6(C)$ gives a section of $g$ over $U$.
\par
\vskip 4pt

{\it However, it is unclear if this section will be extended to a holomorphic section of $g$. The following idea to replace $Y$ further is due to Yasunari Nagai.}
\par
\vskip 4pt

Similarly, one has the moduli space $\overline{X} := M((0, [2H], -4))$ of semi-stable sheaves with Mukai vector $(0, [2H], -4)$ on $S$ and its symplectic resolution $X$ \cite{LS}. For the same reason as before, we have a fibration $\overline{f} : \overline{X} \longrightarrow \mathbf P^5 \simeq \vert 2H \vert$, and the induced fibration $f : X \longrightarrow \mathbf P^5 \simeq \vert 2H \vert$. The fiber over $[C] \in U$ is ${\rm Pic}^0(C)$ for both $\overline{f}$ and $f$. As $g$ admits a section over $U$, we have $f^{-1}(U) \simeq g^{-1}(U)$ over $U$. Thus, $X$ is birational to $Y$, whence, to $\mathcal M_4$. So, $X$ is also deformation equivalent to $\mathcal M_4$ again by [Hu]. 

In order to complete the proof, we need the following slightly 
technical lemma. Note that it is unclear if the structure sheaf of a curve is stable unless the curve is irreducible and reduced. This is pointed out by the referee. 

\begin{lemma} \label{lemma:curve}
For each $C$ in $\vert 2H \vert$, the structure sheaf $\mathcal O_C$ is stable with respect to the polarization $H$. 
\end{lemma}

\begin{proof} As ${\rm Pic}\, S = \mathbf Z H$, an element $C$ of $\vert 2H \vert$ falls into one of the following types (i), (ii), (iii): (i) an irreducible and reduced curve; (ii) a reduced curve of the form $C_1 + C_2$, where $C_1, C_2 \in \vert H\vert$; (iii) a non-reduced curve of the form $2C_1$, where $C_1 \in \vert H \vert$. We put $C_2 = C_1$ for a curve of type (iii).

Let $F$ be a proper subsheaf of $\mathcal O_C$. It suffices to show that $p(F, m) < p(\mathcal O_C, m)$, where $p(*, m)$ is the reduced Hilbert polynomial of a coherent sheaf $*$ with respect to the polarization $H$ (See \cite{HL}, 
pages 10-11). In what follows, we denote 
$F \otimes \mathcal O_{S}(mH)$ by $F(m)$ and so on. When $C$ is of type (i), the result just follows from the exact sequence 
$0 \rightarrow F(m) \rightarrow \mathcal O_C(m) \rightarrow \mathcal O_C/F 
\rightarrow 
0$, where $\mathcal O_C/F$ is a non-zero skyscraper sheaf. 

Let us consider the case where $C$ is of type (ii) or (iii). 

Let $I_{C_i}$ be the ideal sheaf of $C_i$ in $\mathcal O_S$. 
Then $I_{C_1}/ I_{C_1}I_{C_2} \simeq 
\mathcal O_{C_2}(-C_1)$ when $C$ is of type (ii) 
and $I_{C_1}/I_{C_1}^2 \simeq \mathcal O_{C_1}(-C_1)$ when 
$C$ is of type (iii). 
We have then the following exact sequence called the decomposition sequence (see e.g. [BHPV] page 62):
$$0 \rightarrow \mathcal O_{C_2}(-C_1) \xrightarrow{\alpha} \mathcal O_C 
\xrightarrow{\beta} \mathcal O_{C_1} \rightarrow 
0\, .$$ 
Recall that $(H^2) = 2$. Then, by the exact sequence 
$$0 \rightarrow \mathcal O_{S}(m-2) \rightarrow \mathcal O_S(m-1) \rightarrow \mathcal O_{C_2}(-C_1)(m) \rightarrow 
0$$ and by the Riemann-Roch theorem on $S$, we have 
$$\chi (\mathcal O_{C_2}(-C_1)(m)) = 
\chi(\mathcal O_S(m-1)) - \chi(\mathcal O_S(m-2)) = (m-1)^2 - (m-2)^2 = 2m-3\, .$$ 
Similarly, $\chi(\mathcal O_{C}(m)) = 4m -4$ and $\chi(\mathcal O_{C_1}(m)) = 2m-1$, and hence $p(\mathcal O_C, m) = m-1$. 

Put $G = {\rm Im}(\beta \vert F)$ and $K = {\rm Ker}(\beta \vert F)$. 

First we consider the case where $G \not= 0$ and $K \not= 0$. In this case, both $\mathcal O_{C_2}(-C_1)/K$ and $\mathcal O_{C_1}/G$ are skyscraper sheaves. 
We have then 
$$\chi(K(m)) = 2m - 3 - \ell(\mathcal O_{C_2}(-C_1)/K)\, ,\, 
\chi(G(m)) = 2m - 1 - \ell(\mathcal O_{C_1}/G)\, ,$$ 
where $\ell(*)$ is the length of a skyscraper sheaf $*$. Thus
$$\chi(F(m)) = \chi(K(m)) + \chi(G(m)) = 4m-4 - 
(\ell(\mathcal O_{C_2}(-C_1)/K) 
+ \ell(\mathcal O_{C_1}/G))\, ,$$
and hence
$$p(F, m) = m-1 - 
\frac{\ell(\mathcal O_{C_2}(-C_1)/K) 
+ \ell(\mathcal O_{C_1}/G)}{4}\,\,.$$
As, $F \not= \mathcal O_C$, we have $\ell(\mathcal O_{C_2}(-C_1)/K) 
+ \ell(\mathcal O_{C_1}/G) > 0$. Therefore 
$p(F, m) <  p(\mathcal O_C, m)$.

Next we consider the case where $G = 0$. In this case, 
$F = K$ and $\mathcal O_{C_2}(-C_1)/K$ is a skyscraper sheaf. 
Thus $\chi(F(m)) = 2m-3 - \ell(\mathcal O_{C_2}(-C_1)/K)$, and hence 
we have $p(F, m) \le m - 3/2$. Therefore $p(F, m) < p(\mathcal O_C, m)$. 

It remains to consider the case where $K = 0$. When $C$ is of type (ii), we can interchange the roles of $C_1$ and $C_2$ to reduce the problem to the case where $K \not= 0$, and we are done. So, we may assume that 
$C$ is of type (iii). As $F \simeq G$, we have 
$$\chi(F(m)) = \chi(G(m)) = 2m -1 - \ell(\mathcal O_{C_1}/G)$$
and hence
$$p(F, m) = m - \frac{1 + \ell(\mathcal O_{C_1}/G)}{2}\, .$$
Thus $p(F, m) < p(\mathcal O_C, m)$ {\it unless} 
$\ell(\mathcal O_{C_1}/G) = 0$ or $1$. Let us show 
$\ell(\mathcal O_{C_1}/G) \not= 0, 1$ by argue by contradiction. 

Assume that $\ell(\mathcal O_{C_1}/G) = 0$. Then $G = \mathcal O_{C_1}$ and 
the decomposition sequence splits as 
$$\mathcal O_C \simeq \mathcal O_{C_1} \oplus \mathcal O_{C_1}(-C_1)\,$$ 
as $\mathcal O_C$-modules. However, this is impossible, because the right hand side is annihilated by the subsheaf $I_{C_1}/I_{C_1}^2$ of $\mathcal O_C$, while the left hand side is not. 

Assume that $\ell(\mathcal O_{C_1}/G) = 1$. Then 
$G$ is the maximal ideal sheaf ${\bf m}_{C_{1}, P}$ of some point $P$ 
in $C_1$. 
Notice the following natural exact sequence:
$$0 \longrightarrow (f)/(f^2) \longrightarrow (x, y)/(f^2) \longrightarrow (x, y)/(f) \longrightarrow 0\, ,$$
where $(x, y)$ is the local coordinates of $S$ at $P$ and $f$ is the defining equation of $C_1$ at $P$. Then, by restricting the decomposition sequence 
to the maximal ideal ${\bf m}_{C, P}$ of $P$ in $C$, we have the following 
exact sequence of $\mathcal O_C$-modules:
$$0 \rightarrow \mathcal O_{C_1}(-C_1) \xrightarrow{\alpha} {\bf m}_{C, P} 
\xrightarrow{\beta} {\mathbf m}_{C_{1}, P} \longrightarrow 0\, .$$
As $\beta \vert F : F \simeq G = {\mathbf m}_{C_1, P}$, we have then 
$${\mathbf m}_{C, P} \simeq {\mathbf m}_{C_1, P} \oplus 
\mathcal O_{C_1}(-C_1)$$
as $\mathcal O_C$-modules. However, for the same reason as before, this is impossible. 

This completes the proof of Lemma \ref{lemma:curve}.  
\end{proof}

Let us return back to the proof of Theorem \ref{theorem:MAX} (2). Let $B$ be the closed subset of $\overline{X}$ that parametrizes the S-equivalence classes of semistable but unstable sheaves. Then, by Lemma \ref{lemma:curve}, the map $C \mapsto \mathcal O_C$ gives a holomorphic section $s : \mathbf P^5 \longrightarrow P \subset \overline{X}$ of $\overline{f}$ that satisfies $s(\mathbf P^5) \cap B = \emptyset$. As $X$ is just a blow-up of $\overline{X}$ along $B$ \cite{LS}, it follows that $s$ gives rise to a holomorphic section of $f : X \longrightarrow \mathbf P^5$. Now, Theorem \ref{theorem:MAX} (2) follows from Theorem \ref{theorem:DMW}.
\par \vskip 4pt

{\bf Acknowledgement.} First of all, I would like to express my best thanks to Doctor Yasunari Nagai for his brilliant idea in the last step of the proof of Theorem \ref{theorem:MAX} (2) and to the referee for very careful reading and pointing out several mathematical errors and ambiguities in the first version. I would like to express my thanks to Professors F. Catanese, Y. Kawamata, G. Tian and E. Viehweg for their invitation to Oberwolfach September 2007. An initial idea of this work was found during preparation of my talk requested by Professor F. Catanese. I express my thanks to Professors J.H. Keum and J.M. Hwang for invitation to KIAS. Main part of this work was completed during my stay there February 25-April 12 2008.

\end{document}